\documentclass[mathpazo]{article}
\usepackage{authblk}
\usepackage{blindtext}
\usepackage{geometry}
\usepackage{lipsum}
\usepackage[utf8]{inputenc}
\usepackage{graphicx}
\usepackage{epstopdf}
\usepackage{algorithmic}
\ifpdf
  \DeclareGraphicsExtensions{.eps,.pdf,.png,.jpg}
\else
  \DeclareGraphicsExtensions{.eps}
\fi
\usepackage[caption=false]{subfig}
\usepackage{amsmath,amssymb,amsfonts}
\usepackage{verbatim}
\usepackage{xcolor}

\usepackage{amssymb}

\usepackage{amssymb}
\usepackage{amsmath,amssymb,amsfonts}
\usepackage[caption=false]{subfig}
\usepackage{xcolor}

\newcommand{\del}{\partial}
\newcommand{\bs}{\boldsymbol}
\newcommand{\rmd}{\mathrm{d}}
\newcommand{\rmi}{\mathrm{i}}

\newcommand{\rmsc}{\mathrm{sc}}

\newcommand{\rmps}{\mathrm{ps}}

\begin{document}



\title{A non-singular boundary element method for interactions between acoustical field sources and structures\footnote{Published in \emph{Adv. Appl. Math. Mech.} \textbf{15} (2023) 831-851, DOI: 10.4208/aamm.OA-2022-0024.}}

\author{Qiang Sun \thanks{qiang.sun@rmit.edu.au}}
\affil{Australian Research Council Centre of Excellence for Nanoscale Biophotonics, School of Science, RMIT University, Melbourne, VIC 3001, Australia}

%
%
%
\date{}

\maketitle

\begin{abstract}
Localized point sources (monopoles) in an acoustical domain are implemented to a three dimensional non-singular Helmholtz boundary element method in the frequency domain. It allows for the straightforward use of higher order surface elements on the boundaries of the problem. It will been shown that the effect of the monopole sources ends up on the right hand side of the resulting matrix system. Some carefully selected examples are studied, such as point sources near and within a concentric spherical core-shell scatterer (with theoretical verification), near a curved focusing surface and near a multi-scale and multi-domain acoustic lens.
\end{abstract}




\section{\label{sec:1} Introduction}

Sound waves are of importance for numerous applications, ranging from speech to annoying or even hazardous noise prevention. These phenomena can be effectively modelled via the Helmholtz equation. In the frequency domain, assuming a $e^{-\rmi \omega t}$ time dependency for all quantities, with $\omega$ the angular frequency and $t$ time, sound wave phenomena can be described by the Helmholtz equation:

%
\begin{equation}
    \nabla^2 \phi +k^2 \phi = 0, \label{eq:Helmholtz}
\end{equation}
where $\phi$ represents the (velocity) potential and $k$ the wavenumber. The same equation is also valid for the pressure. Eq.~(\ref{eq:Helmholtz}) can be efficiently solved using a boundary element method (BEM). Using BEM for acoustics is still an active area of research \cite{Sumbatyan2021, Zhang2020, Liu2019}.

Often, the effect of the presence of one or a set of (localized) volume sources is important. For example for simulating traffic noise, airplane noise \cite{LugtenJASA2018} or other localized sources of sound \cite{ShengJASA2020}. In those instances, there is no need to model an acoustic source exactly, but it can be simplified as a point source, here represented by an acoustic monopole. To implement these monopoles into numerical methods, due care has to be taken to deal with the physical singularities associated with those point monopoles. 

In this work, we introduce a non-singular boundary element method for acoustics including monopoles in which the mathematically ``artificial'' singularities in the integrands originating from the Green's function are fully removed before any numerical procedures or calculations. In the conventional ways, the mathematical singularities from the Green's function are treated numerically by a local change of variables in the evaluation of the surface integrals~\cite{Telles1987}, which accompanies with additional coding efforts. Some attempts to remove the singularities before the numerical evaluations of surface integrals introduced new unknowns for the tangential derivatives on the surface which leads to unnecessary computations~\cite{Liu1999}. Another idea to deal with the singularities needs to add additional artificial parameters on a dummy ``nearby'' boundary which accuracy strongly depends on the choice of that dummy boundary~\cite{Zhang1999}. Other attempts to treat the mathematical singularities from the Green's function in boundary element methods include density interpolation methods~\cite{PrezArancibia2019, Faria2021} and plane-wave singularity subtraction techniques~\cite{PrezArancibia2018}.

In our non-singular boundary element method (BEM), we show that the physical point monopoles can easily be incorporated in the boundary element method formulations which ends up on the right hand side of the resulting matrix system. Such an analytical treatment is more advantageous than domain methods (volume methods) which need very careful mesh strategies to represent the area near the point source~\cite{Kalayeh2015}. Also, boundary element methods are particularly suited for problems in infinite domains: the radiation conditions at infinity are automatically satisfied and no mesh is needed anywhere else than on the boundaries. Moreover, the non-singular boundary element method presented here takes the implementation a step forward to fully remove the mathematical singular behaviours from the Green's function before any numerical procedures without introducing any unnecessary new unknowns or dummy ``nearby'' boundaries. As such, the regularized integrands allow the easy implementation of high order surface elements to improve the calculation efficiency and accuracy~\cite{Sun2015}. Meanwhile, the non-singular BEM introduced here has great potential to deal with the near singularity issues or the boundary layer effect which is a main issue for the numerical process to evaluate strong or weak singularities of surface integrals by a local change of variables~\cite{Sun2014, Sun2015, Sun2015_Stokes, Sun2016}.


The outline of the paper is as follows: in Sec.~\ref{sec:IncMonopole} the mathematical implementation of the inclusion of a monopole into a fully desingularized boundary element framework is described. The constructed numerical framework is then validated in Sec.~\ref{sec:Validation} for a multi-domain concentric core shell configuration with an analytical solution given in \ref{sec:AppA}. From a physical point of view,  interesting phenomena can arise with the variation of the dimensionless parameter $ka$, where $a$ is a typical dimension of the scattered object. Two examples on focused soundwaves are shown to demonstrate the robustness and effectiveness of our method based on two different physical phenomena: one example with reflection as demonstrated in Sec.~\ref{sec:bowl} and the other with a multi-scale multi-domain acoustic lens as presented in Sec.~\ref{sec:AcousticLens}. Conclusions are given in Sec.~\ref{sec:conclusion}.

\section{\label{sec:theory} Including a monopole in a fully desingularized BEM}\label{sec:IncMonopole}

\begin{figure}[t]
\centering
\includegraphics[width=0.6\textwidth]{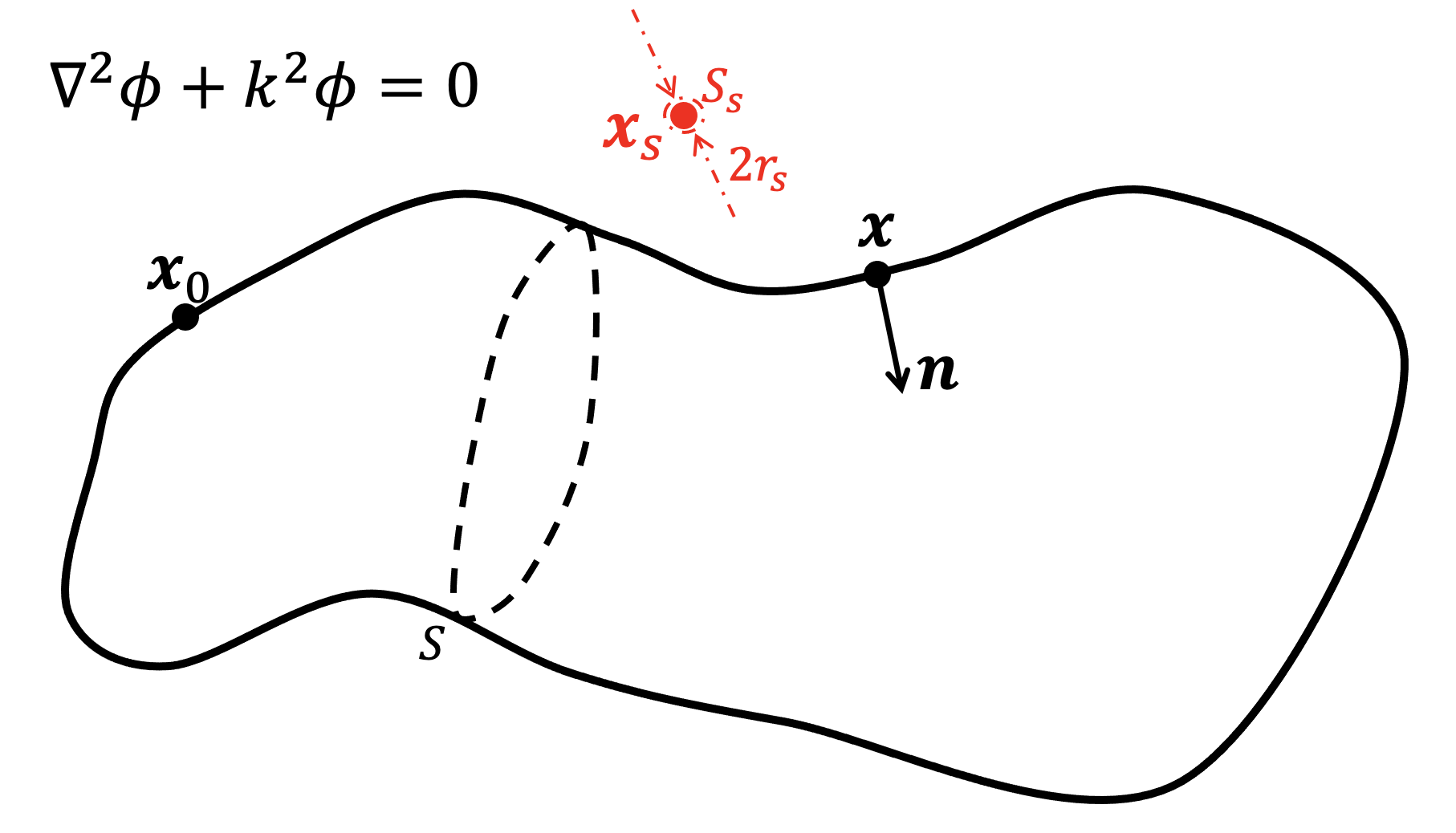} 
\caption{Sketch of the parameters used in the boundary element method with a monopole (located at $\boldsymbol x_S$) included in the domain. The integrals must now be performed on the surface of the scatterer  bounded by the surface $S$ and on a small sphere with radius $r_s$ around the monopole with surface $S_S$. In the theory the limit of $r_s \rightarrow 0$ is taken. A mathematical singularity arises in the conventional boundary element method when $\boldsymbol x$ approaches $\boldsymbol x_0$.}  \label{Fig:0sketch_BEM}
\end{figure}

The main goal of this work is to solve the Helmholtz equation, particularly in the presence of volume point sources, as sketched in Fig.~\ref{Fig:0sketch_BEM}. Since the Helmholtz equation is elliptic, in the standard boundary element framework the Helmholtz equation is rewritten as a surface integral on all surfaces $S$ enclosing the computation domain in terms of the Green's function or kernel $G_k(\boldsymbol x, \boldsymbol x_0) \equiv e^{\rmi k{|\boldsymbol x - \boldsymbol x_0|}}/{|\boldsymbol x - \boldsymbol x_0|}$ between the computation point $\boldsymbol x$ and the observation point $\boldsymbol x_0$ as well as its normal derivative $H_k(\boldsymbol x, \boldsymbol x_0) = \boldsymbol n \cdot \nabla G_k(\boldsymbol x, \boldsymbol x_0)$ with $\bs{n}$ being the unit normal vector pointing out of the computational domain:  
\begin{equation} \label{eq:BEM_original}
    \begin{aligned}
c(\boldsymbol x_0) \phi(\boldsymbol x_0) + \int_{S} \left[ \phi(\boldsymbol x)  H_k(\bs{x},\bs{x}_0) - \frac{\del \phi(\boldsymbol x) } {\del n}  G_k(\bs{x},\bs{x}_0) \right]  \rmd S(\bs{x}) = -I_s.
\end{aligned}
\end{equation}
In Eq.~(\ref{eq:BEM_original}), the constant $c(\boldsymbol x_0)$ is the solid angle at $\boldsymbol x_0$, and $\del \phi(\bs{x}) / \del n \equiv \boldsymbol n \cdot \nabla\phi(\bs{x})$. Eq.~(\ref{eq:BEM_original}) will be applied to all $N$ nodes resulting in a $N \times N$ matrix system. As displayed in Fig.~\ref{Fig:0sketch_BEM}, we have drawn an artificial surface $S_s$ around the monopole source and $I_s$ is defined as: $I_s=\int_{S_s}  [\phi_s  H_k(\bs{x},\bs{x}_0) - (\del \phi_s / \del n)  G_k(\bs{x},\bs{x}_0) ] \rmd S(\bs{x})$, with $\phi_s$ the potential near the monopole and $\boldsymbol x$ is now located on $S_s$, while $\boldsymbol x_0$ is still situated on $S$. Thus the singular behavior of the integral $I_s$ is now due to the singular behavior of $\phi$ instead of $H_k$ and $G_k$. The term $I_s$ will be determined later. 

There are now two kinds of singularities in the problem. The first kind is of  mathematical origin and is associated with the Green's function when the integration point $\boldsymbol x$ approaches the point $\boldsymbol x_0$ on the surface $S$ of the problem. The second (physical) singularity is due to the monopole in the domain. We will show that both singularities can be relatively easily dealt with in a fully non-singular boundary element framework.

Let us deal first with the mathematical singularity caused by the Green's function. There are several ways to desingularize this equation by subtracting a suitable similar solution that has the same singular behavior on the surface $S$, for example 
\begin{equation} \label{eq:BEM_Laplace}
    \begin{aligned}
c(\boldsymbol x_0) \psi(\boldsymbol x_0) + \int_{S} \left[  \psi(\boldsymbol x)  H_0(\bs{x},\bs{x}_0) - \frac{\del \psi(\boldsymbol x) } {\del n}  G_0(\bs{x},\bs{x}_0) \right]  \rmd S(\bs{x}) = 0
\end{aligned}
\end{equation}
with $G_0(\boldsymbol x, \boldsymbol x_0) = 1/|\boldsymbol x - \boldsymbol x_0|$ and $H_0(\boldsymbol x, \boldsymbol x_0) = -(\boldsymbol x -\boldsymbol x_0) \cdot \boldsymbol n /|\boldsymbol x - \boldsymbol x_0|^3$ the Green's function for a 
Laplace problem and its normal derivative. Since Eq.~(\ref{eq:BEM_Laplace}) represents a Laplace equation, $\nabla^2 \psi = 0$, a constant and a linear function are a solution. If $\psi$ is chosen to be the following combination of a constant and linear function 
\begin{equation} \label{eq:psiFunction}
\psi(\boldsymbol x) = \phi(\boldsymbol x_0) + \boldsymbol n(\boldsymbol x_0) \cdot (\boldsymbol x - \boldsymbol x_0) \frac{\del \phi(\boldsymbol x_0)}{\del n},
\end{equation}
then the function $\psi$ will approach $\phi(\boldsymbol x_0)$ when $\boldsymbol x$ approaches $\boldsymbol x_0$. Also, when the surface $S$ is smooth then $\del \psi(\boldsymbol x_0)/\del n = \boldsymbol n(\boldsymbol x_0) \cdot \boldsymbol n(\boldsymbol x) \del \phi(\boldsymbol x_0) / \del n \rightarrow  \del \phi(\boldsymbol x_0) / \del n$. When Eq.~(\ref{eq:BEM_Laplace}) is subtracted from Eq.~(\ref{eq:BEM_original}), a totally non-singular boundary element framework emerges as: 
\begin{equation} \label{eq:BEM_design}
    \begin{aligned}
4\pi \phi(\boldsymbol x_0) + \int_{S}  \left[\phi(\boldsymbol x)  H_k(\bs{x},\bs{x}_0) -\psi(\boldsymbol x)  H_0(\bs{x},\bs{x}_0)\right] \rmd S(\bs{x}) \\
= \int_{S} \left[\frac{\del \phi(\boldsymbol x) } {\del n}  G_k(\bs{x},\bs{x}_0) - \frac{\del \psi(\boldsymbol x) } {\del n}  G_0(\bs{x},\bs{x}_0) \right] \rmd S(\bs{x}) - I_s.
\end{aligned}
\end{equation}
The factor $4 \pi \phi(\boldsymbol x_0)$ that appears in Eq.~(\ref{eq:BEM_design}) originates from the fact that the integral at the surface at infinity can no longer be ignored for infinite domains, for instance the example demonstrated in Fig.~\ref{Fig:0sketch_BEM}. The term $\phi(\boldsymbol x_0)$ in Eq.~(\ref{eq:psiFunction}) is responsible for this term. Note that the solid angle has disappeared from Eq.~(\ref{eq:BEM_design}). Also, it is worth noticing that the singular behavior of $G_k$ is in fact the same as $G_0$ since $G_k \equiv G_0 + \Delta G$ where $\Delta G \equiv (e^{\rmi k{|\boldsymbol x - \boldsymbol x_0|}}-1)/{|\boldsymbol x - \boldsymbol x_0|}$ is finite when $\bs{x}\rightarrow \bs{x}_0$, and the same analysis can be applied to $H_{k}$. The approach shown in Eq.~(\ref{eq:BEM_design}) was tested thoroughly in the past for Laplace problems and Stokes flow problems~\cite{Klaseboer2012, Sun2013, Sun2014, Sun2015_Stokes}, Helmholtz equations~\cite{Sun2015,2112.08693}, molecular and colloidal electrostatics~\cite{Sun2016}, electromagnetic scattering ~\cite{Klaseboer2017,Sun2017,KlaseboerAO2017,Sun2020a,Sun2020b} and interactions between light and matter~\cite{Sun2021,Sun2022}. The desingularization using Eq.~(\ref{eq:BEM_Laplace}) is by no means the only way to desingularize the classical boundary element method, many alternative methods can be found~\cite{Klaseboer2012, Sun2015}.

Since the desingularized boundary element method with high order elements, such as quadratic curved elements, is so accurate, a further advantage is that it is very unlikely to hit an internal resonance frequency and get a so-called spurious solution (for a more elaborate discussion on the fictitious frequency problem see~\cite{KlaseboerEABE2019}). Thus, there is no need to apply Burton-Miller type of schemes~\cite{BurtonMiller}. However, if needed, the non-singular Burton-Miller BEM is also available~\cite{2112.08693}.

The next task is to determine the integral $I_s$ over the surface $S_s$ surrounding the small volume with the monopole source. The velocity potential of a monopole with constant complex strength $Q$ is $\phi_s= Q e^{\rmi kr_s}/(4 \pi r_s)$, where $Q$ represents the volume of fluid displaced by the source~\cite{Kinsler2000}, then $\partial \phi_s/\partial n=\boldsymbol n \cdot \nabla \phi_s = -Q \rmd  [e^{\rmi kr_s}/(4 \pi r_s)]/ \rmd r_s$. In the neighbourhood of $S_s$, this potential will totally overpower any other potentials and thus: 
\begin{equation}
\label{eq:Is1}
    \begin{aligned}
    I_s &= \lim_{r_s \rightarrow 0}\int_{S_s}  \left[\phi_s  H_k(\bs{x},\bs{x}_0) - \frac{\del \phi_s } {\del n}  G_k(\bs{x},\bs{x}_0)\right]  \rmd S(\bs{x})  \\
   &\approx \lim_{r_s \rightarrow 0}\int_{S_s} Q  \frac{e^{\rmi kr_s}}{4 \pi r_s}   H_k(\bs{x_s},\bs{x}_0) \rmd S(\bs{x}) +  \lim_{r_s \rightarrow 0}\int_{S_s} Q \frac{\rmd } {\rmd r_s}\left[ \frac{e^{\rmi kr_s}}{4 \pi r_s}\right]   G_k(\bs{x_s},\bs{x}_0) \rmd S(\bs{x}).
    \end{aligned}
\end{equation}
Here we have used that $H_k(\boldsymbol x, \boldsymbol x_0)$ and $G_k(\boldsymbol x, \boldsymbol x_0)$ are not singular, since $\boldsymbol x$ is situated on the surface $S$ and $\boldsymbol x_0$ on $S_s$. For a very small sphere, $H(\boldsymbol x, \boldsymbol x_0)$ can thus be replaced by $H_k(\boldsymbol x_s, \boldsymbol x_0)$ and similar for $G_k$. The first integral with $H_k$ behaves as $1/r_s \cdot r_s^2\rightarrow 0$, since $\rmd S\sim r_s^2$, and thus vanishes. The second integral becomes:
\begin{equation}
    \begin{aligned}
       \lim_{r_s \rightarrow 0}\int_{S_{s}}  Q \frac{\rmd} {\rmd r_s}\left[ \frac{e^{\rmi kr_s}}{4 \pi r_s} \right]   \rmd S(\bs{x})&=\lim_{r_s \rightarrow 0}\int_{S_{s}} Q\frac{e^{\rmi kr_s}}{4 \pi r_s^2} (\rmi kr_s -1)  \rmd S(\bs{x})
      =-Q 
    \end{aligned}
\end{equation}
since $e^{\rmi kr_s} (\rmi kr_s -1)/r_s^2 \sim -1/r_s^2$ and $\rmd S=4\pi r_s^2$. Thus:
\begin{equation}\label{eq:Is2}
    \begin{aligned}
    I_s = \lim_{r_s \rightarrow 0}\int_{S_{s}} \left[ \phi_s  H_k(\bs{x},\bs{x}_0) - \frac{\del \phi_s } {\del n}  G_k(\bs{x},\bs{x}_0) \right]  \rmd S(\bs{x})  = -Q \frac{e^{\rmi k|\boldsymbol x_s - \boldsymbol x_0|}}{|\boldsymbol x_s - \boldsymbol x_0|}.
    \end{aligned}
\end{equation}
The above proof can easily be extended to multiple acoustic sources with different strengths placed at different locations in the domain. For example, for $M$ monopoles one gets:
\begin{equation} \label{eq:BEM_designSum}
    \begin{aligned}
4\pi \phi(\boldsymbol x_0) + \int_{S}  \left[ \phi(\boldsymbol x)  H_k(\bs{x},\bs{x}_0) -\psi(\boldsymbol x)  H_0(\bs{x},\bs{x}_0) \right ] \rmd S(\bs{x}) \\
 =\int_{S} \left [\frac{\del \phi(\boldsymbol x) } {\del n}  G_k(\bs{x},\bs{x}_0) - \frac{\del \psi(\boldsymbol x) } {\del n}  G_0(\bs{x},\bs{x}_0)  \right ] \rmd S(\bs{x}) + \sum_{i=1}^M Q_i \frac{e^{\rmi k|\boldsymbol x_{s,i} - \boldsymbol x_0|}}{|\boldsymbol x_{s,i} - \boldsymbol x_0|}.
\end{aligned}
\end{equation}
Eq.~(\ref{eq:BEM_designSum}) is the boundary regularised integral equation formulation including monopoles for acoustics used in this work. It is robust, efficient and accurate for solving multi-scale multi-domain acoustic problems. In particular, it is user-friendly and easy to be implemented together with high order elements since no treatment is needed to deal with singularities. To demonstrate this, a few examples are presented in Sec.~\ref{sec:Validation} to Sec.~\ref{sec:AcousticLens}. In all those simulations, 6-noded quadratic triangular elements~\cite{Sun2016} are used to represent the surfaces and potentials. 

Note that the above described method (including the desingularization) can directly be applied to problems involving the Laplace equation ($k=0$) as well. Once $\phi$ and $\partial \phi/ \partial n$ are known on the surface after Eq.~(\ref{eq:BEM_designSum}) is solved, through postprocessing the values of $\phi$ anywhere in the domain can be obtained. This can also be done in a desingularized manner, for more details see \cite{Sun2015}.

In our soundwave simulations, we use three different kinds of boundary conditions. In the `hard surface' scattering, it is assumed that the normal velocity on the surface of the scatterer is zero ($\partial \phi/ \partial n = 0$). For simulations with two or more acoustic media, the boundary conditions across an interface are the continuity of the normal velocity (thus $\partial \phi/\partial n$ must be continuous across the interface) and the pressure (which assumes that $\rho \phi$ is continuous across the interface, with $\rho$ the density), originating from the relationship $ p= \rmi \omega \rho \phi$ (with $\omega$ the angular frequency of the system, see \cite{KirkupBook}. Finally, for a bubble, it is assumed that the pressure at the bubble surface is zero, and thus $\phi=0$ at the bubble surface.

\section{\label{sec:result} Validation with a core shell model} \label{sec:Validation}

\begin{figure}[!ht]
\centering
\subfloat[]{\includegraphics[width=0.4\textwidth]{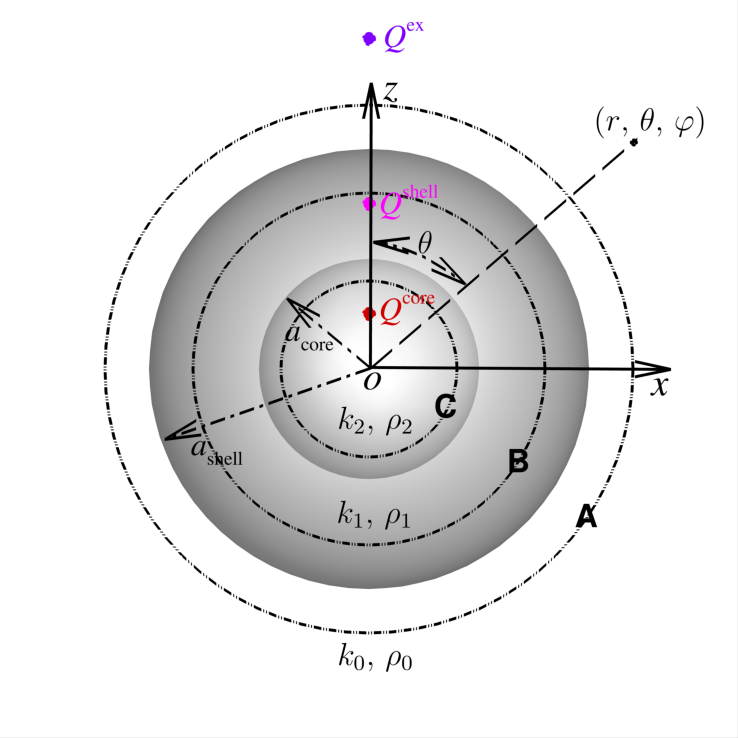}} \quad
\subfloat[]{\includegraphics[width=0.4\textwidth]{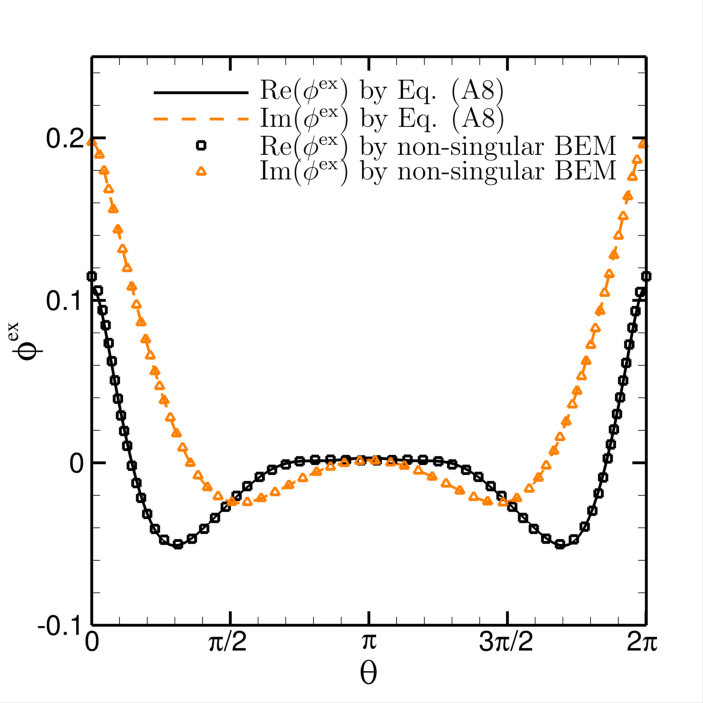}} \\
\subfloat[]{\includegraphics[width=0.4\textwidth]{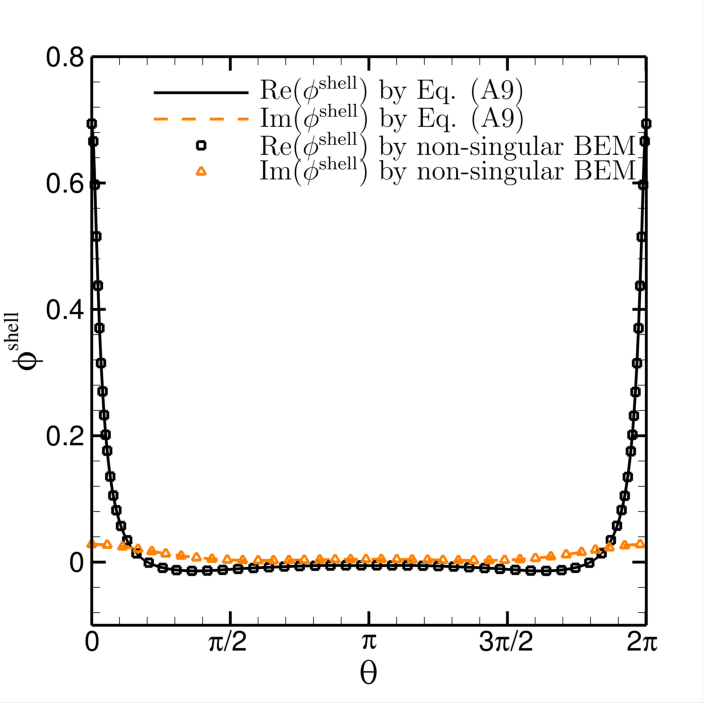}} \quad
\subfloat[]{\includegraphics[width=0.4\textwidth]{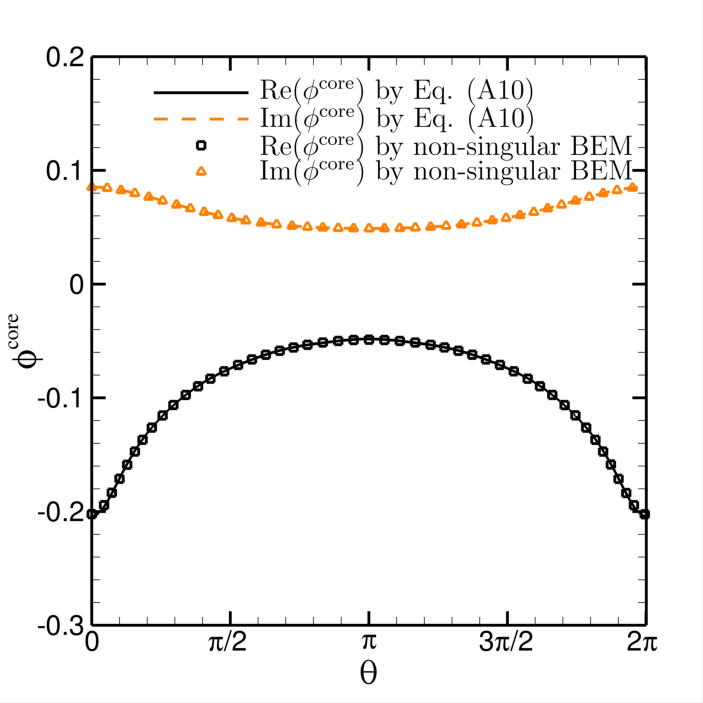}}
\caption{Excellent agreement is shown between the analytical results obtained from Eqs.~(\ref{eq:anaex}) to (\ref{eq:anacore}) and the numerical results obtained by using the boundary element method for the potential in the external domain (b), in the shell (c) and in the core (d) when a concentric spherical core-shelled scatterer is excited by three acoustic monopoles. The analytical results are indicated by lines and the numerical results by symbols.} \label{Fig:0anacompare}
\end{figure}

To validate our non-singular boundary element method with domain monopoles included, we consider a concentric spherical core-shell scatterer in an external medium as shown in Fig.~\ref{Fig:0anacompare}(a). There is one monopole in the scatterer core, $Q^{\text{core}}$ one in the scatterer shell, $Q^{\text{shell}}$ and one in the surrounding medium (the external domain), $Q^{\text{ex}}$ to drive the acoustic field simultaneously. These three monopoles are located along the axis of symmetry of the spheres, the $z$-axis of the Cartesian coordinate system with its origin locating at the centre of the core-shell spherical scatterer. In App.~\ref{sec:AppA}, the analytical solution of such a problem is obtained. 

The physical and geometrical parameters are: for the radii ratio, $a_{\text{shell}}/a_{\text{core}}=2$; with the wavenumber for the external domain $k_0$ satisfying $k_0a_{\text{core}}=1$, the ratio of wavenumbers of the inner and shell are $k_1/k_0=1.5$ and $k_2/k_0=0.8+\rmi 0.6$, respectively. Note, that $k_2$ has been chosen to be complex (i.e it will caused damping). The source strengths are chosen to be (note that we have deliberately chosen two of them out of phase and one with an imaginary number to validate our method) $Q^{\text{ex}}=0.8+\rmi 0.6$ at $r_s^{\text{ex}}/a_{\text{core}}=3$, $Q^{\text{shell}}=1.0$ at $r_s^{\text{shell}}/a_{\text{core}}=1.5$, $Q^{\text{core}}=-1.0$ at $r_s^{\text{core}}/a_{\text{core}}=0.5$. Finally the density ratios of the core and shell with respect to the outer medium are $\rho_1/\rho_0=5$ and $\rho_2/\rho_0=2$, respectively. The comparison with the theory has been done on three circular tracks, centered at the origin, one in the external domain, one in the shell and one in the core, all in the $xz$ plane, as displayed in Fig.~\ref{Fig:0anacompare} (a) with dashed lines, indicated with {\bf{A}}, {\bf{B}} and {\bf{C}}. The radii of the these three tracks are $r_{\textbf{A}}/a_{\text{core}}=2.4$, $r_{\textbf{B}}/a_{\text{core}}=1.6$ and $r_{\textbf{C}}/a_{\text{core}}=0.8$. For the case with the parameters mentioned above, the number of terms of the sum in Eqs.~(\ref{eq:anaex}) to (\ref{eq:anacore}) is calculated to be $N=6$. The values of each of the wavenumbers $k$ are specially selected to make sure that $N$ for each monopole is comparable. As shown in Fig.~\ref{Fig:0anacompare}(b)-(d), excellent agreement is obtained between the analytical results from Eqs.~(\ref{eq:anaex}) to (\ref{eq:anacore}) and the numerical results when 642 nodes connected by 320 quadratic elements are used for the shell and the same number for the core surface. The numerical and theoretical solutions are virtually overlapping. Also, the convergence and accuracy tests on our non-singular BEM have been performed. As shown in Fig.~\ref{Fig:app_sketch} (b), when the observation point is located at $r=1.2a_{\text{shell}}$ along the $x$-axis from the centre of the concentric spherical core-shell scatterer, the relative error of the potential between the numerical results by the non-singular BEM and the analytical solution in Eq.~(\ref{eq:anaex}) reduces from $10^{-4}$ for 1284 nodes (including the nodes on the shell and the core) to about $3\times 10^{-7}$ for 15684 nodes. The computational time increases from a few seconds to 900 s.

The numerical implementation of this multi-domain (3-domain) boundary element problem consists of a few steps. First, we discretise the two surfaces, say $N_{\text{shell}}$ nodes on the interface between the shell and the external domain and $N_{\text{core}}$ on the surface in between the shell and the core particle. For the external domain we then get $N_{\text{shell}}$ equations from Eq.~(\ref{eq:BEM_design}), one for each node $\boldsymbol x_0$. It appears most convenient to express all variables in terms of the shell potential. With the boundary conditions we thus get $\phi^{\text{ex}}=\phi^{\text{shell}} \rho_1/\rho_0$ and $\partial \phi^{\text{ex}}/\partial n = \partial \phi^{\text{shell}}/\partial n$. For the shell domain a similar boundary element method will be applied for each node on the surfaces that enclose the computation domain which are the shell surface and the core surface, as such now $N_{\text{shell}} + N_{\text{core}}$ in total. Also, $k_0$ is now replaced by $k_1$, the term with $4 \pi$ is no longer present since it is an internal domain, and the source term must be replaced. Finally, for the internal domain another $N_{\text{core}}$ equations will be obtained with wavenumber $k_2$. This will result in a matrix system of size $2(N_{\text{shell}} +  N_{\text{core}}) \times 2(N_{\text{shell}} +  N_{\text{core}})$ which can now be solved to get the unknowns $\phi^{\text{shell}}$ and its normal derivative $\partial \phi^{\text{shell}}/\partial n$ on both the shell surface and core surface.

\begin{figure}[t]
\centering
\subfloat[]{\includegraphics[width=0.4\textwidth]{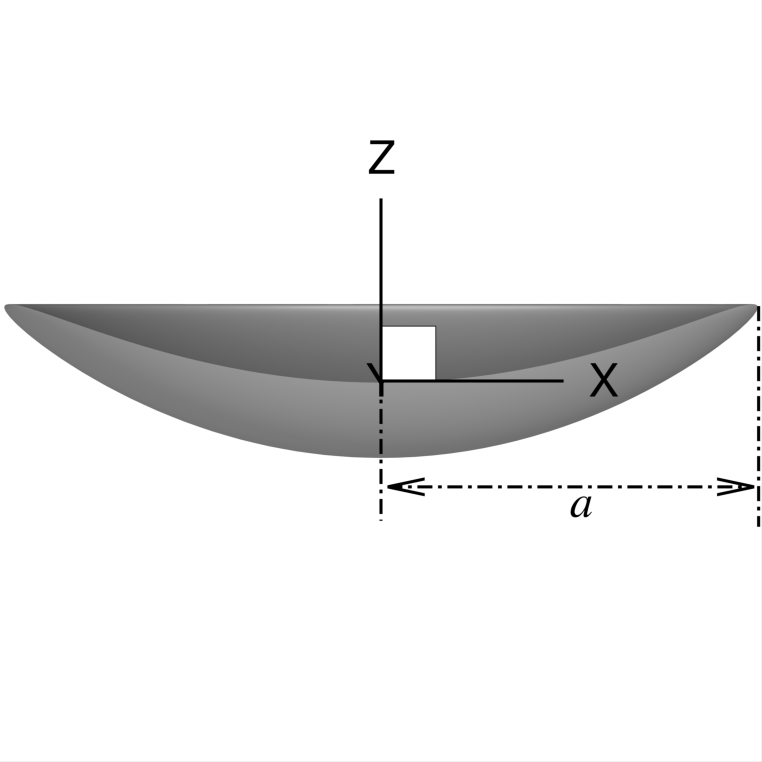}} \qquad
\subfloat[]{\includegraphics[width=0.4\textwidth]{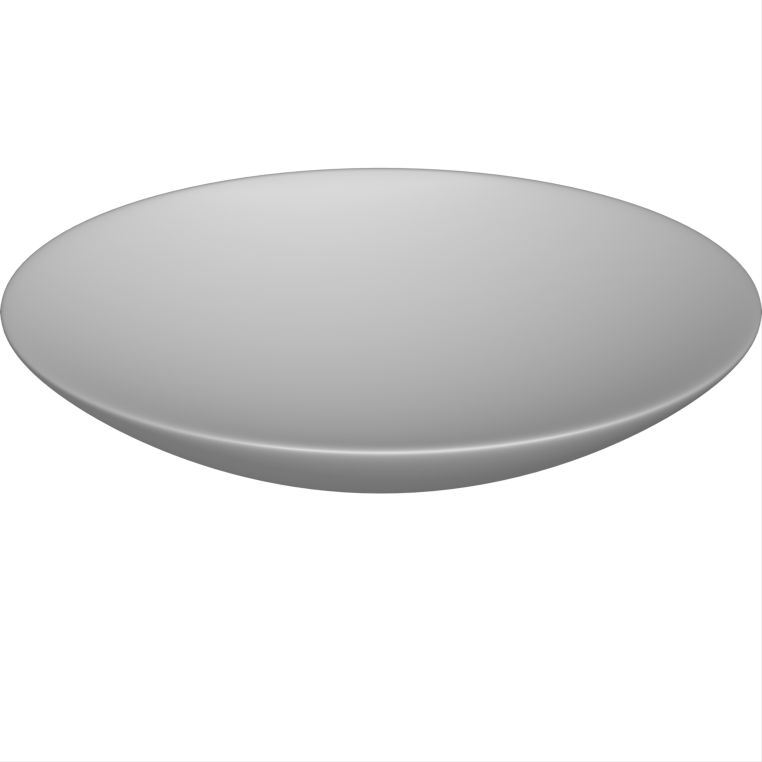}}
\caption{Bowl shaped acoustic reflector given by $(x,y,z) = (a \sin{\theta}\cos{\varphi}, \, a \sin{\theta}\sin{\varphi}, \, 0.1a (\cos{\theta}-1)+0.3\sin^2{\theta})$ with $\theta$ the polar angle and $\varphi$ the azimuthal angle: (a) Cross sectional view revealing the top and bottom surfaces of the reflector and (b) Full 3D view.} \label{Fig:bowlshape}
\end{figure}
\begin{figure}[!ht]
\centering
\subfloat[]{\includegraphics[width=0.32\textwidth]{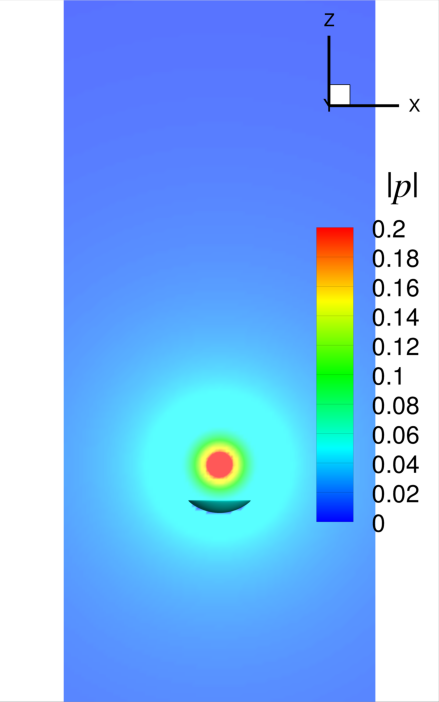}}
\subfloat[]{\includegraphics[width=0.32\textwidth]{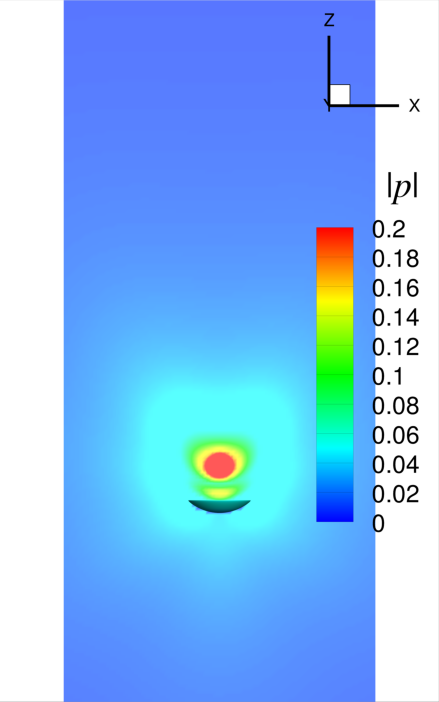}}
\subfloat[]{\includegraphics[width=0.32\textwidth]{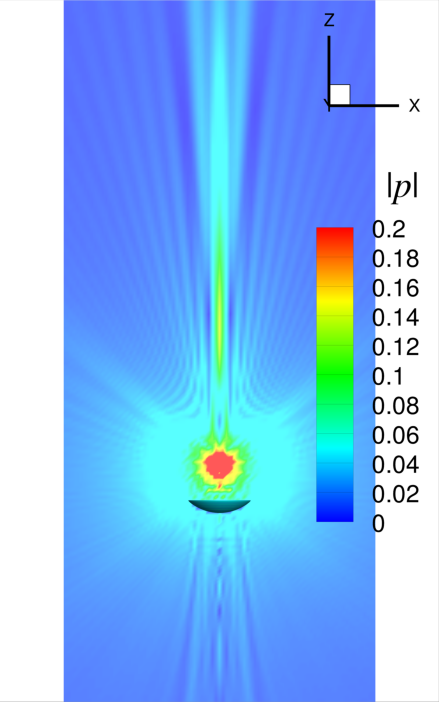}} \\ 
\subfloat[]{\includegraphics[width=0.32\textwidth]{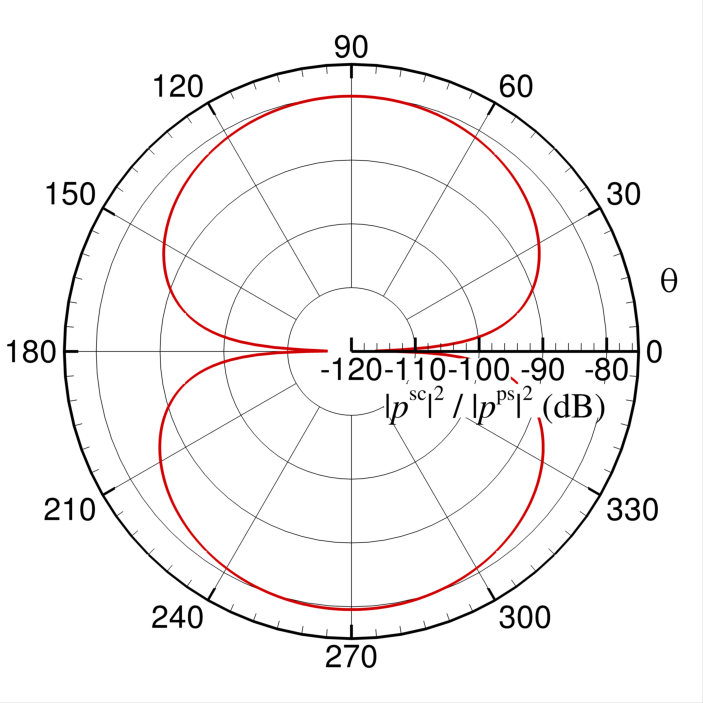}}
\subfloat[]{\includegraphics[width=0.32\textwidth]{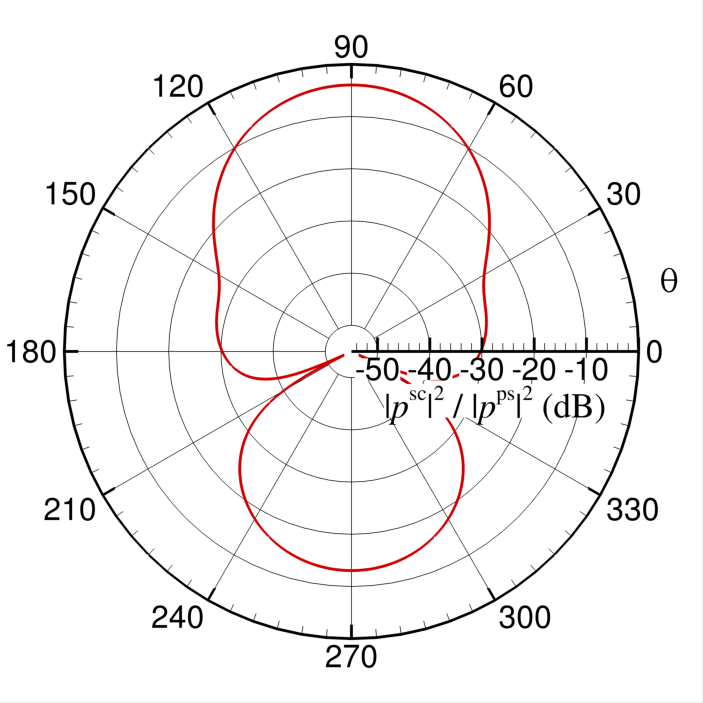}}
\subfloat[]{\includegraphics[width=0.32\textwidth]{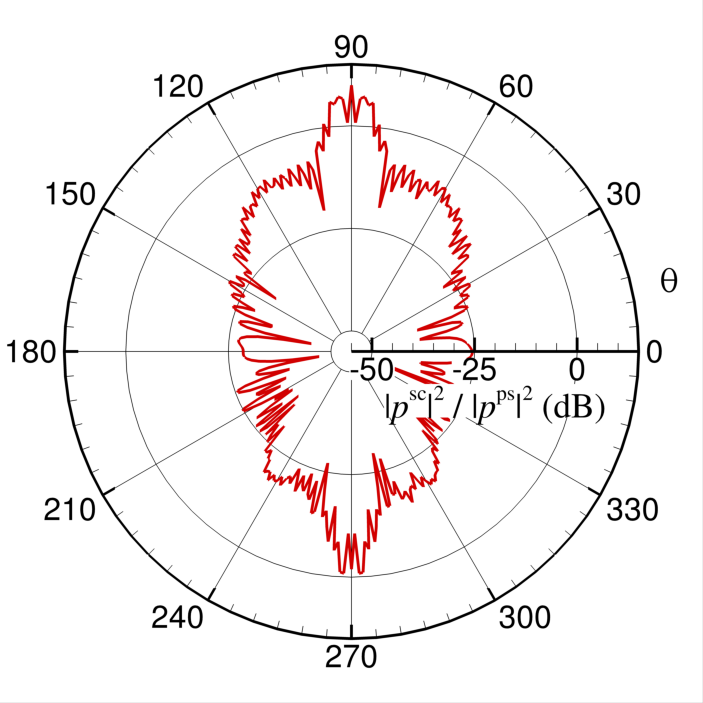}}
\caption{Simulation of the bowl shaped object from Fig~\ref{Fig:bowlshape} with a monopole placed at the focal point of the bowl. The horizontal scale of the bowl is $2a$. For $ka=0.01$, the wave essentially does not `see' the bowl and the radar plot is (d) essentially becomes that of a dipole. For $ka=5$ (b) and (e) a marked reflection can be observed. For $ka=100$, (c) and (f) a very strong upward reflecting beam can be observed. In (a)-(c), the pressure is normalized by the pressure of the monopole at $r=|Q\rho_0\omega|/(4\pi)$ in the external domain. The radar plots (d), (e) and (f) were taken at $r=100 a$.} \label{Fig:1bowlka}
\end{figure}

\section{One or more monopoles near a bowl shaped object} \label{sec:bowl}

Now that it is clear that the constructed numerical framework is reliable and accurate, we will show some further examples of interesting situations with a monopole interacting with a scatterer for which no analytical solutions are available. 

We start with a monopole source placed next to the concave surface of a rigid bowl shaped object in a medium with density of $\rho_0$ and the scattering from this object is simulated. It is interesting to investigate the effect of the wavelength on the scattered pattern with respect to the size of the object (thus changing the variable $ka$). The bowl shape is given by $(x,y,z) = (a \sin{\theta}\cos{\varphi}, \, a \sin{\theta}\sin{\varphi}, \, 0.1a (\cos{\theta}-1)+0.3\sin^2{\theta})$ with $\theta$ the polar angle and $\varphi$ the azimuthal angle. The monopole is located along the axis of symmetry of the bowl at $(0, \, 0, \, 1.35a)$, as presented in Fig.~\ref{Fig:bowlshape}. At the surface of the rigid bowl, the normal velocity of the fluid is zero. The non-singular boundary element method including one monopole was used to calculate the pressure distribution in the surrounding area of the bowl and the far field scattered pattern when the wavelength of the acoustic wave generated by the monopole is much larger than, in the same order of, and much smaller than the size of the bowl. In these calculations, the bowl surface is represented by 2880 quadratic triangular elements connected by 5762 nodes. When the acoustic wavelength, $\lambda$ with $k \lambda = 2\pi$, from the monopole is much larger than the size of the bowl as $ka = 0.01 $, the acoustic wave is barely affected by the presence of the bowl. As such, the distribution of pressure, $p$, in the area surrounding the bowl is very similar to that of a single monopole in free space (external domain), as shown in Fig.~\ref{Fig:1bowlka}(a), and in the far field, the scattered field, $p^{\rmsc}$, is a dipole field as shown in Fig.~\ref{Fig:1bowlka}(d). Here, $p^{\rmsc}$ is calculated by $p=p^{\rmsc}+p^{\rmps}$ with $p^{\rmps}$ the pressure of the monopole in free space. When the wavelength of the acoustic field generated by the monopole is of the same order as the size of the bowl with $ka=5$, the effect of the bowl on the acoustic field becomes obvious as the field is reflected and the field in front of the bowl is elongated as shown in both Figs.~\ref{Fig:1bowlka}(b) and~\ref{Fig:1bowlka}(e). When the wavelength is much smaller relative to the bowl size with $ka=100$, ray phenomena are becoming dominant. Since the monopole is located at the focal point of the concave bowl, a focused beam along the axis of symmetry appears, and the reflected acoustic field from the bowl is much stronger relative to the diffraction field at the back of the bowl, as shown in Figs.~\ref{Fig:1bowlka}(c) and~\ref{Fig:1bowlka}(f).

\begin{figure}[t]
\centering
\subfloat[]{\includegraphics[width=0.32\textwidth]{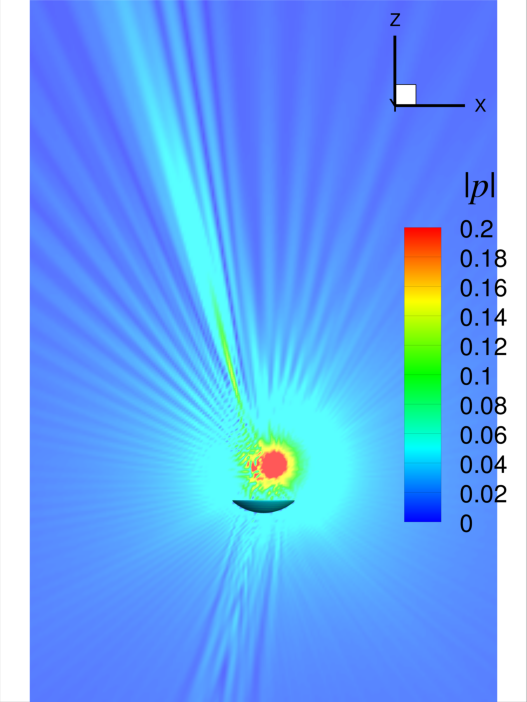}}
\subfloat[]{\includegraphics[width=0.32\textwidth]{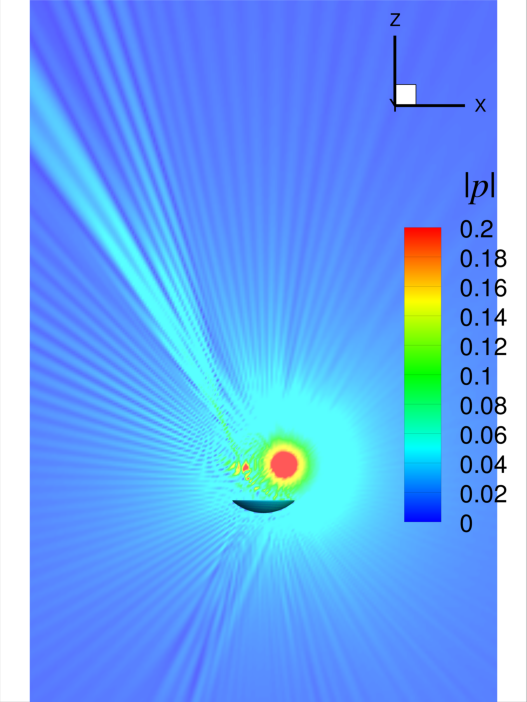}}
\subfloat[]{\includegraphics[width=0.32\textwidth]{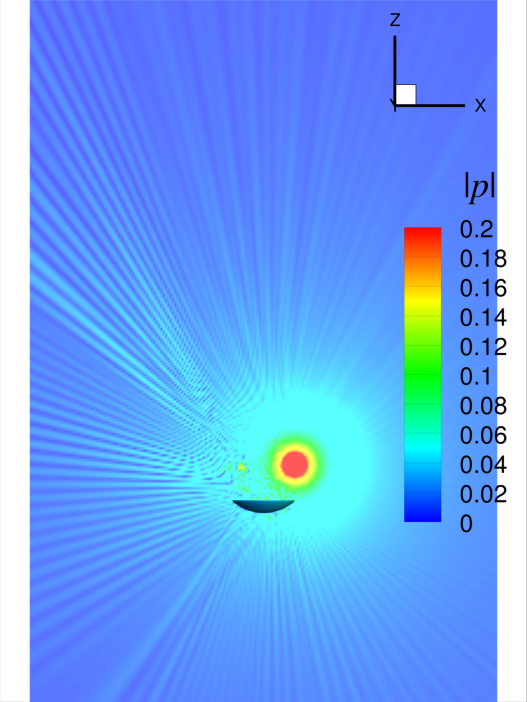}}
\caption{The effect of changing the location of the monopole with respect to the bowl shaped object ($ka=100$). The monopole is placed $0.33 a$ (a), $0.67a$ (b) and $1.0 a$ (c) in the x-direction instead of on the axis of symmetry. the pressure is normalized by the pressure of the monopole at $r=|Q\rho_0\omega|/(4\pi)$ in the external domain.} \label{Fig:2bowlka100off}
\end{figure}

When the acoustic monopole is placed away from the symmetry axis of the blow with $ka=100$, the reflection patterns were affected accordingly, as shown in Fig~\ref{Fig:2bowlka100off}. When the monopole location is shifted horizontally to $(x_s, \, y_s, \, z_s) = (0.33a, \, 0, \, 1.35a)$, the reflected focused beam is rotated counterclockwise around the focal point of the concave bowl by about 13.75$^{\text{o}}$ (Fig.~\ref{Fig:2bowlka100off}a). When the monopole is moved further away from the symmetric axis of the bowl to $(x_s, \, y_s, \, z_s) = (0.67a, \, 0, \, 1.35a)$, the focused beam rotates more in the counterclockwise direction (Fig.~\ref{Fig:2bowlka100off}b). If the monopole is moved to the edge of the bowl at $(x_s, \, y_s, \, z_s) = (1.0a, \, 0, \, 1.35a)$, the focused beam almost vanishes (Fig.~\ref{Fig:2bowlka100off}c). 

\begin{figure}[t]
\centering
\subfloat[]{\includegraphics[width=0.32\textwidth]{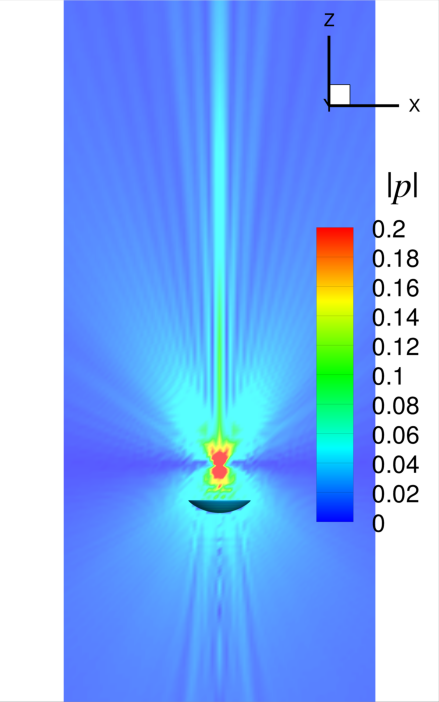}}
\subfloat[]{\includegraphics[width=0.32\textwidth]{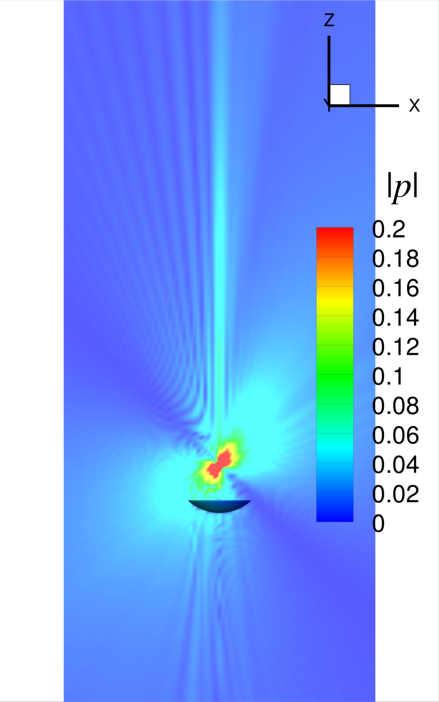}}
\subfloat[]{\includegraphics[width=0.32\textwidth]{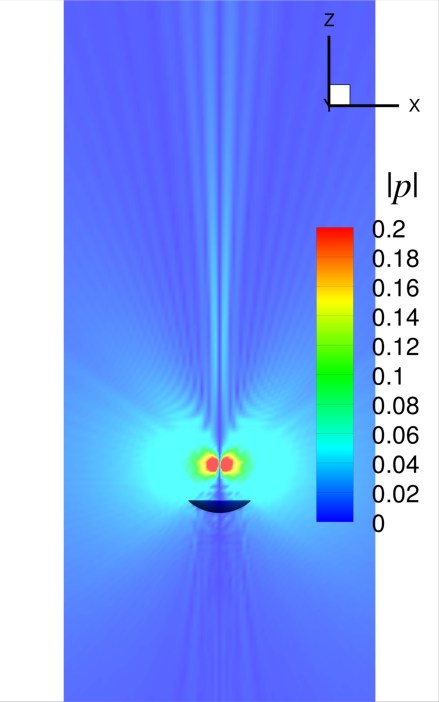}}
\caption{Dipole at the focal point of the bowl shaped object. (a) dipole aligned with z-axis, (b) at a 45 degrees angle and (c) at a 90 degrees angle. In all three cases $ka=100$. the pressure is normalized by the pressure of the dipole at $r=|Q\rho_0\omega|/(4\pi)$ in the external domain.} \label{Fig:3bowldipole}
\end{figure}

A dipole can be constructed by two monopoles with opposite strength (fully out of phase) located very close to each other (here taken to be at a distance $10^{-5} a$). If this dipole is located at the focal point of the concave bowl, the acoustic pressure patterns are shown in Fig.~\ref{Fig:3bowldipole}. Since acoustic waves are longitudinal, when the dipole axis is aligned with the symmetry axis of the bowl, a strong focussed beam is formed from the reflected field of the bowl along its symmetric axis and the field is fully symmetric with respect to the bowl symmetry axis, as shown in Fig.~\ref{Fig:3bowldipole}(a). When the dipole axis is rotated $\pi/4$ clockwise with respect to the $z$-axis, the symmetry of the acoustic field is broken, as shown in Fig.~\ref{Fig:3bowldipole}(b). If the dipole axis is perpendicular to the symmetric axis of the rigid bowl, the symmetry of the field with respect to that axis is recovered with a low pressure field along that axis and two relative high pressure fields next to the axis due to the focus effect from the reflection field of the concave surface of the rigid bowl, as shown in Fig.~\ref{Fig:3bowldipole}(c).

\section{Multi-scale multi-domain acoustic lens} \label{sec:AcousticLens}

\begin{figure}[t]
\centering
\subfloat[]{\includegraphics[width=0.4\textwidth]{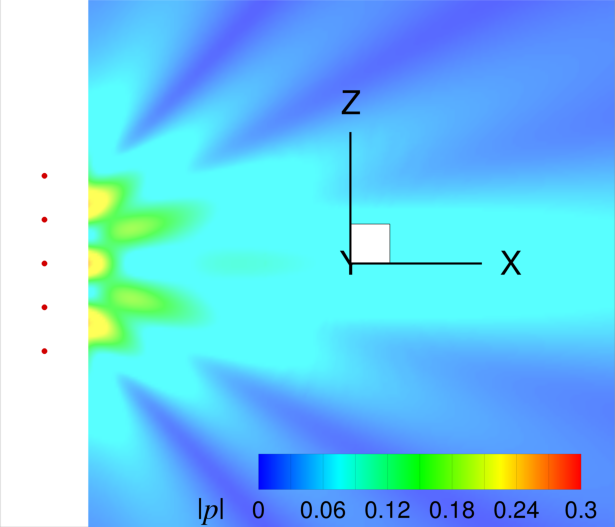}}\quad \quad
\subfloat[no core]{\includegraphics[width=0.4\textwidth]{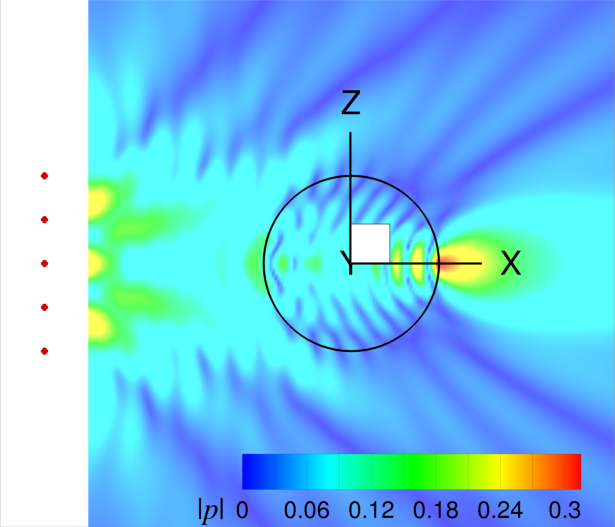}} \\ 
\subfloat[rigid core]{\includegraphics[width=0.4\textwidth]{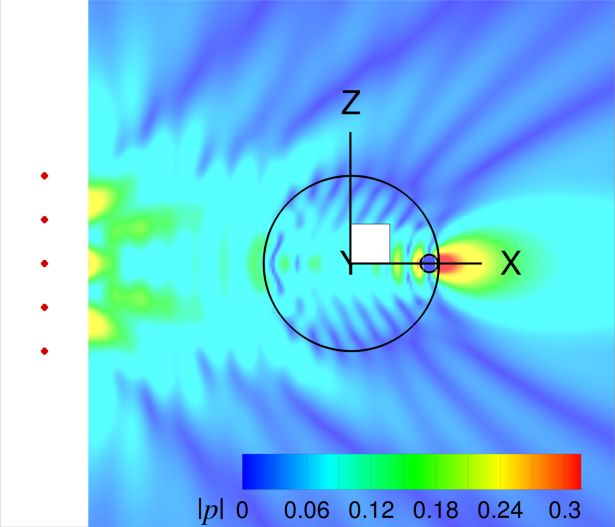}} \quad \quad
\subfloat[air bubble core]{\includegraphics[width=0.4\textwidth]{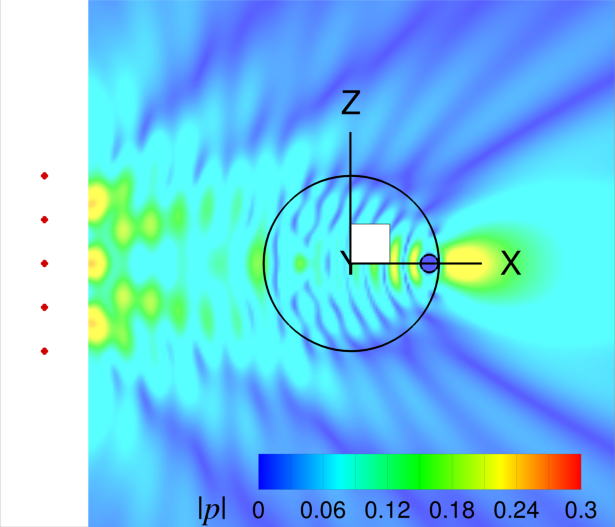}}
\caption{Pressure amplitudes for 5 monopoles in water with (a) no drop, (b) an oil drop, (c) an oil drop with a tiny rigid sphere at the back and (d) an oil drop with a tiny air bubble at the back. The parameter are: $k_0a=10$, $\rho_{\text{drop}}/\rho_0$ = 0.8, $k_{\text{drop}}/k_0$ = 1.333. The five monopoles are located at $x = -3.5a_{\text{drop}}$, $y=0$, and $z=-a_{\text{drop}}$, $-0.5a_{\text{drop}}$, $0$, $0.5a_{\text{drop}}$ and $a_{\text{drop}}$,  respectively (indicated with red dots). Maxima and minima at about half the wavelength appear due to interference of scattered and incoming waves appear in (b), (c) and (d), while a more or less uniform amplitude is observed in the area around $x$=0 in (a).} \label{Fig:4coreshell}
\end{figure}

The next interesting example is an acoustic lens, comprising of a spherical oil drop of radius $a_{\text{drop}} =$ 2.39 mm submerged in water, with or without a tiny rigid sphere or bubble included. The ensemble is being irradiated by five acoustic monopoles with frequency 1 MHz. The monopoles are aligned perpendicular to the drop as indicated in Fig.\ref{Fig:4coreshell}(b) by the red dots. They are placed $h^{\rmps} =$ 1.19 mm apart. The distance between the line of monopoles and the center of the oil drop is 8.4 mm (or $3.5 a_{\text{drop}}$). The reference pressure amplitude obtained without the presence of the oil drop is shown in Fig.\ref{Fig:4coreshell}(a). Next, the oil droplet is placed at the origin of the coordinate system in Fig.\ref{Fig:4coreshell}(b). The presence of this oil droplet clearly creates a focal point at the back the droplet ($\max{|p|}=0.322$ at $x=0.980 a_{\text{drop}}$). The wavenumber in the external domain (water) is noted as $k_0$ with $k_0 a_{\text{drop}} = 10$ when the speed of sound in water is 1500 m/s. The density ratio of the oil drop to water is $\rho_{\text{drop}}/\rho_{0} = 0.8$ and the ratio of the wavenumber is $k_{drop} / k_0 = 1.33$. The effects of surface tension have been ignored in the simulations. Both the surfaces of the oil drop and the small defect inside it are represented by 720 quadratic triangular elements connected by 1442 nodes. The pressure is normalized by the pressure of a single monopole at $r=|Q\rho_0\omega|/(4\pi)$ in the external domain.

\begin{figure}[t]
\centering
\subfloat[]{\includegraphics[width=0.32\textwidth]{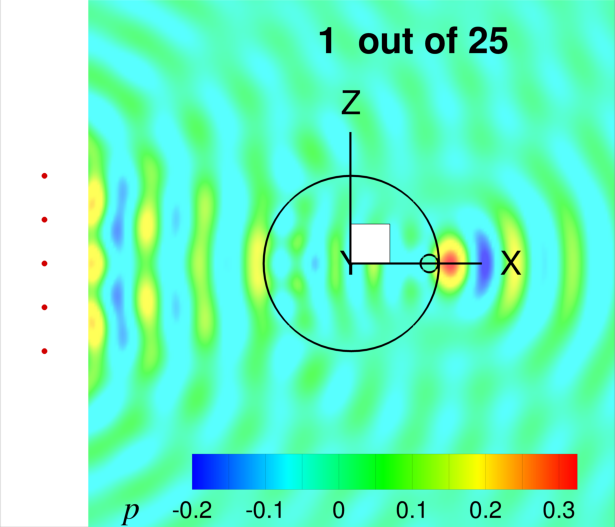}}
\subfloat[]{\includegraphics[width=0.32\textwidth]{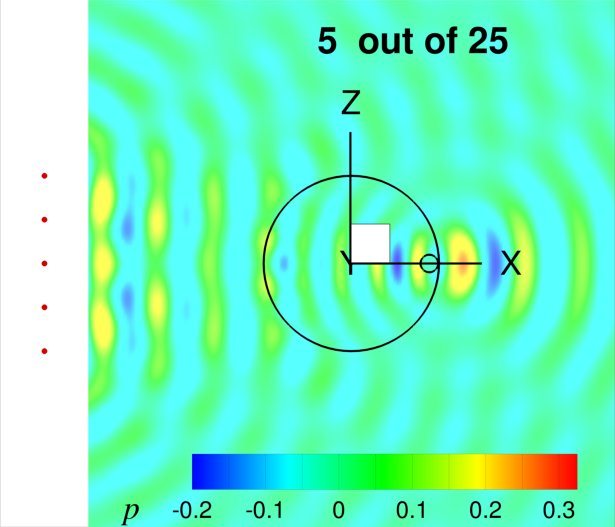}}
\subfloat[]{\includegraphics[width=0.32\textwidth]{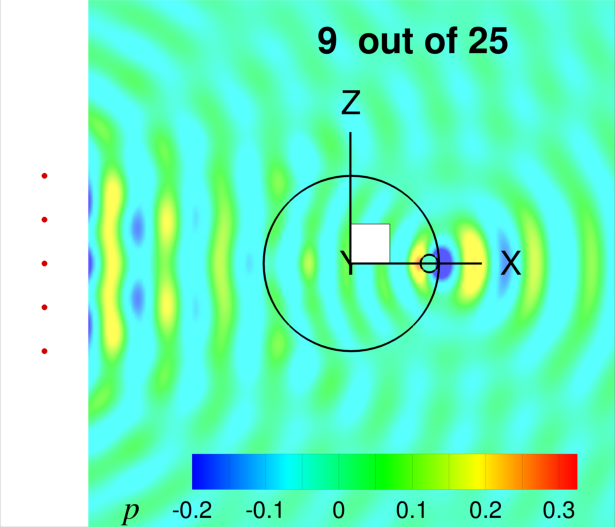}} \\
\subfloat[]{\includegraphics[width=0.32\textwidth]{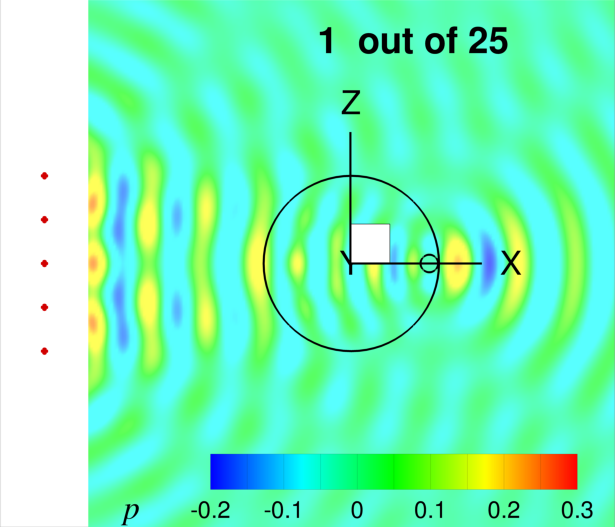}}
\subfloat[]{\includegraphics[width=0.32\textwidth]{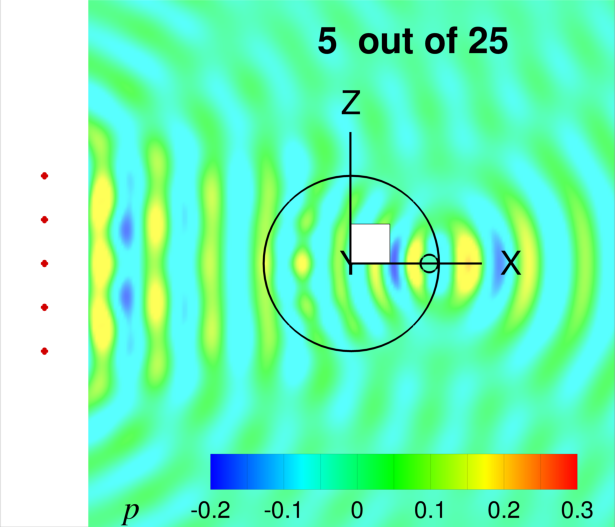}}
\subfloat[]{\includegraphics[width=0.32\textwidth]{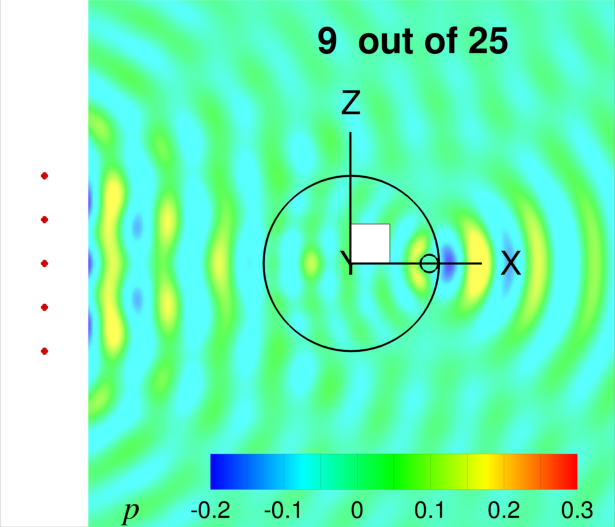}}
\caption{As Fig.\ref{Fig:4coreshell}, but instantaneous pressure profile at selected times: (a)-(c), rigid core; (d)-(f) bubble core. The pressure is normalized by the pressure of a single monopole at $r=|Q\rho_0\omega|/(4\pi)$ in the external domain. Note the recombining spherical waves from the monopoles (the red dots) into a plane wave region in front of the drop. The interplay between reflected and transmitted waves results in a focal point just behind the drop. Full animations are available as supplementary material. } \label{Fig:5time}
\end{figure}

It is instructive to investigate if this acoustic `lens' can be affected by the presence of some defect inside the oil drop, a multi-scale multi-domain acoustic lens. For example the inclusion of a tiny particle at the back of the oil droplet changes noticeably the scattering pattern as clearly illustrated in Fig.\ref{Fig:4coreshell}(c), the focusing effect is stronger ($\max{|p|}=0.398$ at $x=0.980 a_{\text{drop}}$) when compared to that of Fig.~\ref{Fig:4coreshell}(b). Here, a rigid particle with radius $a_{\text{core}}/a_{\text{drop}} = 0.1$ is added with the minimum gap between the inner core and the oil shell $h_{\text{gap}}/a_{\text{drop}} = 0.01$. If an air bubble with the same size of that tiny particle is placed at the same location, the focal point is moved further away from the drop as presented in Fig.~\ref{Fig:4coreshell}(d) and its amplitude is clearly decreased ($\max{|p|}=0.244$ at $x=1.221 a_{\text{drop}}$). Thus an object smaller than the wavelength of sound can still have an observable influence on the formation of a focal point. 

Snapshots of the time-harmonic pressure distributions are shown in Fig.~\ref{Fig:5time} for both the oil drop with the rigid sphere and with the bubble. It is clear that the combined wave field of the five sources combined to a pseudo plane wave according to the Huygens principle before it hits the oil drop. The full animations are available as supplementary material.

\section{\label{sec:conclusion} Conclusions}

Often the effect of a sound source can be regarded as a localized point source. In this work, we demonstrated a non-singular boundary element method for the acoustics with point monopoles. As shown in Sec.~\ref{sec:IncMonopole}, one advantage of our method is that the monopoles can be relatively straightforward included in the non-singular boundary element framework for the Helmholtz equation. Also, as the solid angle and the singularities in the integrands have been fully removed before any numerical procedures or calculations, the non-singular boundary element method introduced here can be implemented with high order surface elements, such as quadratic triangular elements used in Secs.~\ref{sec:Validation} to~\ref{sec:AcousticLens}, straightforwardly. Another advantage of our non-singular boundary element method is that it does not introduce any unnecessary new unknowns or dummy ``nearby'' boundaries. Meanwhile, the non-singular BEM demonstrated here has great potential to deal with the near singularity issues or the boundary layer effect which is a big challenge for the numerical treatments of strong or weak singularities by a local change of variables in the evaluation of surface integrals. Our future work is to combine the non-singular boundary element method presented here with fast methods, such as fast multipole method~\cite{Liu2014}, adaptive cross approximation~\cite{Kurz2007} and hierarchical off-diagonal low-rank matrix method~\cite{Huang2017}, to improve the computational efficiency for practical applications. We believe that our non-singular boundary element framework can make a good contribution to the development and usage of boundary element methods for solving practical multi-scale multi-domain acoustic problems with point volume monopoles. 

\section*{Acknowledgments}
Q.S. thanks Dr Evert Klaseboer from A-STAR Institute of High Performance Computing (Singapore) for the discussions and suggestions on this work. This work was partially supported by the Australian Research Council (ARC) through Grants DE150100169, FT160100357 and CE140100003.

\appendix

\section{Core-shell scatterer with monopoles; theory} \label{sec:AppA}

\begin{figure}[t]
\centering
\subfloat[]{\includegraphics[width=0.45\textwidth]{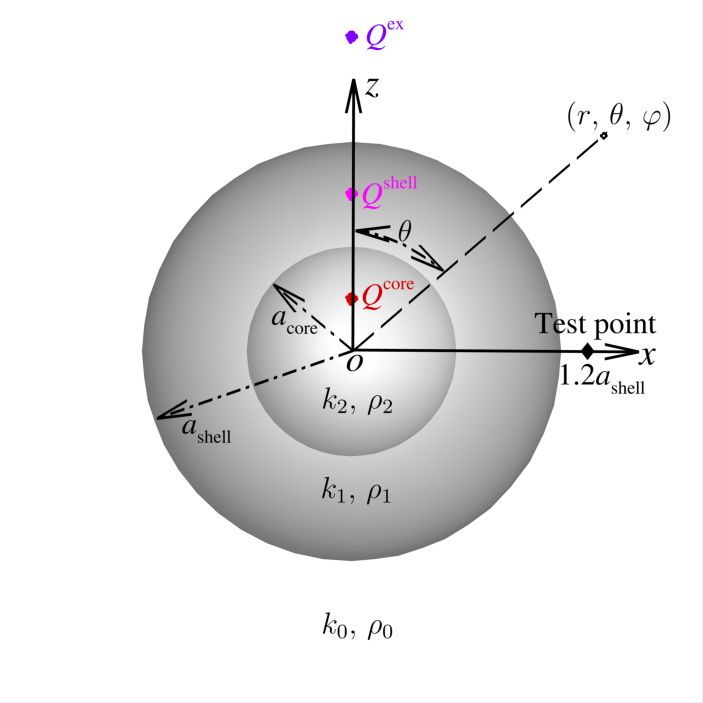}}\qquad
\subfloat[]{\includegraphics[width=0.45\textwidth]{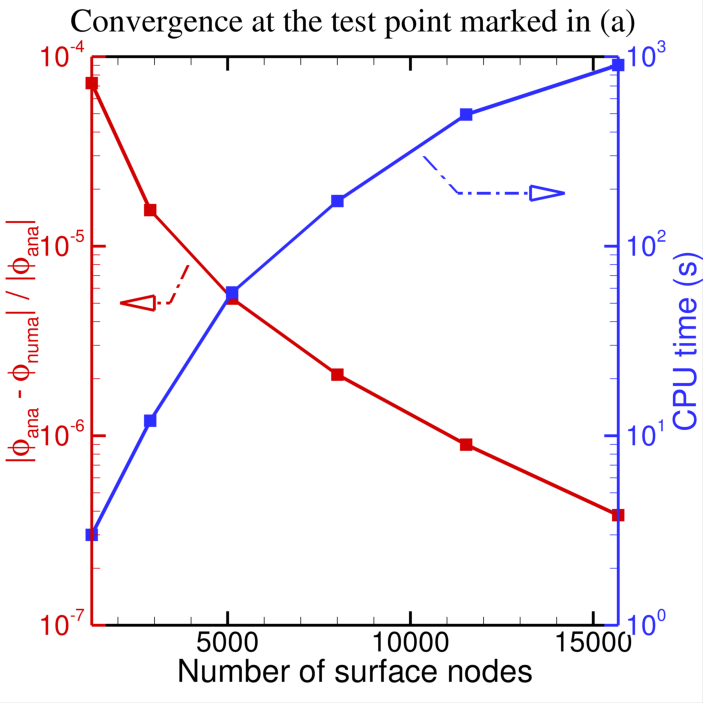}}
\caption{(a) Sketch of the physical problem when a concentric spherical core-shell scatterer is irradiated by three acoustic monopoles all located along the z-axis (the axis of symmetry). A Cartesian coordinate system and the corresponding spherical coordinate system are set with their origin at the centre of the spherical core. (b) Convergence test for the observation point at $r=1.2a_{\text{shell}}$ along the $x$-axis from the centre of the concentric spherical core-shell scatterer. The relative error between the numerical results by the non-singular BEM and the analytical solution in Eq.~(\ref{eq:anaex}) reduces from $10^{-4}$ for 1284 nodes (including the nodes on the shell and the core) to about $3\times 10^{-7}$ for 15684 nodes. The computational time increases from a few seconds to 900 s.}  \label{Fig:app_sketch}
\end{figure}

In this appendix, the analytical solution of the field driven by the acoustic monopoles located along the axis of symmetry of a concentric core-shell spherical scatterer is derived. As shown in Fig.~\ref{Fig:app_sketch} (a), the centre of a core-shell spherical scatterer is set at the origin of a Cartesian coordinate system. The axis of symmetry of this spherical scatterer is chosen along the $z$-axis where the acoustic monopoles are also situated. To solve the acoustic field of such a system, it is most convenient to employ a spherical coordinate system which origin is at the center of the spherical scatterer and polar angle is measured from the $z$-axis. 

Suppose that an acoustic point source (monopole) is situated at $\bs{x}_s \equiv (0, \, 0, \, r_s)$, (with $r_s=|\boldsymbol x_s|$). The expression of the acoustic field at $\bs{x} \equiv (r, \, \theta, \, \varphi)$, $\phi^{\rmps}(\bs{x}) \equiv \phi^{\rmps}(r,\theta)$, generated by such a monopole does not depend on the azimuthal angle $\varphi$ and can be represented asymptotically in terms of free spherical multipolar waves~\cite{Burke1995} as 
\begin{align}\label{eq:anaGF_spherical}
    \phi^{\rmps}(\bs{x}) \equiv \phi^{\rmps}(r,\theta)  = & \frac{Q}{4\pi}\frac{\exp{(\rmi k|\bs{x}-\bs{x}_s|)}}{|\bs{x}-\bs{x}_s|} \nonumber \\
    = & \frac{Q}{4\pi}\rmi k \sum^{N}_{n=0} (2n + 1) h_{n} (k r_{>}) j_n(k r_{<}) P_{n}(\cos{\theta})
\end{align}
where $Q$ is the strength of the acoustic monopole, $\bs{x}$ is the position vector, $r_{>} \equiv \max(r, \, r_s)$, $r_{<} \equiv \min(r, \, r_s)$, $k$ is the domain wavenumber, $h_{n} (k r_{>})$ is the spherical Hankel function of the first kind, $j_n(k r_{<})$ is the spherical Bessel function of the first kind, and $N$ is number of terms in the sum (determining the truncated error). Inspired by the Mie scattering theory, see \cite{BohrenHuffman}, we propose that $N$ can be estimated using the following formula to balance the computational efficiency and stability with the truncated error (lower than 10$^{-3}$):
\begin{align}
    N = |k| + 4|k|^{\frac{1}{3}}+1  & \qquad \text{if } |\bs{x}_s| \ge \frac{1}{2}; \nonumber \\
    N = |(k|\bs{x}_s|)| + 4|(k|\bs{x}_s|)|^{\frac{1}{3}} +1 & \qquad \text{if } |\bs{x}_s| < \frac{1}{2}
\end{align}  
in which the expression on the right hand side is rounded off to the nearest integer.

Also, consider a scalar wave equation for function $\phi$ with wavenumber $k$:
\begin{align}\label{eq:anawaveeq}
    \nabla^2 \phi + k^2 \phi = 0.
\end{align}
Eq.~(\ref{eq:anawaveeq}) is variable separable in spherical coordinates, and its elementary solutions are in the form
\begin{subequations}\label{eq:anawaveeq_sol}
    \begin{align}
    \phi_{(l,n)} = &\cos{(l \varphi)} P_{n}^{l}(\cos{\theta}) z_{n}(kr), \\
    \phi_{(l,n)} = &\sin{(l \varphi)} P_{n}^{l}(\cos{\theta}) z_{n}(kr)
    \end{align}
\end{subequations}
where $l$ and $n$ are integers ($n \ge l \ge 0$), $P_{n}^{l}(\cos{\theta})$ is an associated Legendre polynomial, and $z_{n}(kr)$ is the spherical Bessel function. The following rules are applied to determine the choice of the function $z_{n}(kr)$. In the bounded domain comprising the origin, $z_{n}(kr) \equiv j_{n}(kr)$, the spherical Bessel function of the first kind, is used since $j_{n}(kr)$ is finite at the origin. In the bounded domain excluding the origin, both $z_{n}(kr) \equiv j_{n}(kr)$ and $z_{n}(kr) \equiv y_n(kr)$, the spherical Bessel functions of the first and second kind, are needed. In the unbounded external domain, the spherical Hankel function $z_{n}(kr) \equiv h_{n}(kr) = j_{n}(kr) + \rmi y_n(kr)$ is used since $\rmi k h_{n}(kr) \sim \rmi^{n} \exp{(\rmi k r)}/r$, representing an outgoing wave. A good reference for all the special functions mentioned above can be found in NIST Digital Library of Mathematical Functions~\cite{NISTDLMF}, and some useful recurrence relations in our derivation are listed below: 
\begin{subequations}\label{eq:anasBf}
    \begin{align}
    z_{n+1} (\zeta) &= \frac{2n+1}{\zeta} z_{n}(\zeta) -z_{n-1}(\zeta) \qquad \qquad n=1,2,...\qquad , \\
    \frac{\rmd z_{n}(\zeta)}{\rmd \zeta} & \equiv z_{n}'(\zeta)  =-z_{n+1}(\zeta) + \frac{n}{\zeta}z_{n}(\zeta) \qquad n=0,1,... \qquad ;
    \end{align}
\end{subequations}
\begin{subequations}\label{eq:anasBf012}
    \begin{align}
    j_0(\zeta) &= \frac{\sin\zeta}{\zeta},  \\
    j_1(\zeta) &= \frac{\sin\zeta}{\zeta^2} - \frac{\cos\zeta}{\zeta}, \\
    j_2(\zeta) &= \left(-\frac{1}{\zeta} + \frac{3}{\zeta^3}\right) \sin\zeta - \frac{3}{\zeta^2} \cos\zeta, \\
    y_0(\zeta) &= -\frac{\cos\zeta}{\zeta},  \\
    y_1(\zeta) &= -\frac{\cos\zeta}{\zeta^2} - \frac{\sin\zeta}{\zeta}, \\
    y_2(\zeta) &= \left(\frac{1}{\zeta} - \frac{3}{\zeta^3}\right) \cos\zeta - \frac{3}{\zeta^2} \sin\zeta;
    \end{align}
\end{subequations}
and
\begin{subequations}
    \begin{align}
        P_{0}(\cos{\theta}) & = 1, \\
        P_{1}(\cos{\theta}) & = \cos{\theta}, \\
        P_{n+2}(\cos{\theta}) & =  \frac{(2n+3)\cos{\theta}}{n+2}P_{n+1}(\cos{\theta}) -\frac{n+1}{n+2}P_{n}(\cos{\theta}).
    \end{align}
\end{subequations}

When comparing the expressions  Eq.~(\ref{eq:anaGF_spherical}) and Eq.~(\ref{eq:anawaveeq_sol}), and considering the discussion below Eq.~(\ref{eq:anaGF_spherical}), we can write for the acoustic field potential, in the external domain when $r \geq a_{\text{shell}}$ 
\begin{align}\label{eq:anaex}
\phi^{\text{ex}}(r,\theta) =& \quad \frac{Q^{\text{ex}}}{4\pi}\rmi k_{0} \sum^{N}_{n=0} (2n + 1) h_{n} (k_{0} r^{\text{ex}}_{>}) j_n(k_{0} r^{\text{ex}}_{<}) P_{n}(\cos{\theta})  \nonumber \\
& +  \sum^{N}_{n=0} C_{n} h_{n}(k_{0} r) P_{n}(\cos{\theta}),
\end{align} 
in the spherical shell when $a_{\text{core}} \leq r \leq a_{\text{shell}} $ 
\begin{align}\label{eq:anashell}
\phi^{\text{shell}}(r,\theta) =& \quad \frac{Q^{\text{shell}}}{4\pi}\rmi k_{1} \sum^{N}_{n=0} (2n + 1) h_{n} (k_{1} r^{\text{shell}}_{>}) j_n(k_{1} r^{\text{shell}}_{<}) P_{n}(\cos{\theta})  \nonumber \\
& +  \sum^{N}_{n=0} \left[D_{n} j_{n}(k_{1} r) + E_{n} y_{n}(k_{1} r)\right] P_{n}(\cos{\theta}),
\end{align} 
and in the spherical core when $r \leq a_{\text{core}}$ 
\begin{align}\label{eq:anacore}
\phi^{\text{core}}(r,\theta) = &\quad \frac{Q^{\text{core}}}{4\pi}\rmi k_{2} \sum^{N}_{n=0} (2n + 1) h_{n} (k_{2} r^{\text{core}}_{>}) j_n(k_{2} r^{\text{core}}_{<}) P_{n}(\cos{\theta})  \nonumber \\
& +  \sum^{N}_{n=0} F_{n} j_{n}(k_{2} r) P_{n}(\cos{\theta}).
\end{align}
In Eqs.~(\ref{eq:anaex}) to (\ref{eq:anacore}), $k_0$, $k_1$ and $k_2$ are the wavenumber in the external domain, the shell and the core, respectively; $r^{\text{ex}}_{>} \equiv \max(r, \, r^{\text{ex}}_s)$, $r^{<} \equiv \min(r, \, r^{\text{ex}}_s)$, $r^{\text{shell}}_{>} \equiv \max(r, \, r^{\text{shell}}_s)$, $r^{<} \equiv \min(r, \, r^{\text{shell}}_s)$, $r^{\text{core}}_{>} \equiv \max(r, \, r^{\text{core}}_s)$, $r^{<} \equiv \min(r, \, r^{\text{core}}_s)$, respectively; and $C_{n} \, D_{n}, \, E_{n}, \, F_{n}$ are the unknown coefficients to be determined by the boundary conditions. 

To calculate the unknown coefficients in Eqs.~(\ref{eq:anaex}) to (\ref{eq:anacore}), the boundary condition of the continuity of the acoustic pressure and that of the continuity of the radial velocity are used on the interface between the core and shell and that between the shell and external domain (see also the discussion on boundary conditions in Section~\ref{sec:theory}). As such, we have
\begin{subequations}\label{eq:anabc}
    \begin{align}
        \rho_0 \phi^{\text{ex}}(a_{\text{shell}},\theta) &= \rho_1 \phi^{\text{shell}}(a_{\text{shell}},\theta), \\
        \rho_1 \phi^{\text{shell}}(a_{\text{core}},\theta) &= \rho_2 \phi^{\text{core}}(a_{\text{core}},\theta), \\
        \frac{\partial \phi^{\text{ex}}(r,\theta)}{\partial r}\bigg|_{r=a_{\text{shell}}}  &= \frac{\partial \phi^{\text{shell}}(r,\theta)}{\partial r}\bigg|_{r=a_{\text{shell}}}, \\
        \frac{\partial \phi^{\text{shell}}(r,\theta)}{\partial r}\bigg|_{r=a_{\text{core}}}  &= \frac{\partial \phi^{\text{core}}(r,\theta)}{\partial r}\bigg|_{r=a_{\text{core}}},
    \end{align}
\end{subequations}
where $\rho_0$, $\rho_1$ and $\rho_2$ are the density of the external domain, the shell and the core, respectively. Introducing Eqs.~(\ref{eq:anaex}) to (\ref{eq:anacore}) into Eq.~(\ref{eq:anabc}) leads to 
\begin{subequations}\label{eq:anabcCDEF}
    \begin{align}
        & \rho_0 C_{n} h_{n}(k_{0} a_{\text{shell}}) - \rho_1[D_{n} j_{n}(k_{1} a_{\text{shell}}) + E_{n} y_{n}(k_{1} a_{\text{shell}})] \nonumber \\
        =& -\rho_0 \frac{Q^{\text{ex}}}{4\pi}\rmi k_{0} (2n + 1) h_{n} (k_{0} r^{\text{ex}}_{s}) j_n(k_{0} a_{\text{shell}}) + \rho_1 \frac{Q^{\text{shell}}}{4\pi}\rmi k_{1} (2n + 1) h_{n} (k_{1} a_{\text{shell}}) j_n(k_{1} r^{\text{shell}}_{s}) \\
        & \rho_1[D_{n} j_{n}(k_{1} a_{\text{core}}) + E_{n} y_{n}(k_{1} a_{\text{core}})] - \rho_2F_{n} j_{n}(k_{2} a_{\text{core}}) \nonumber \\
        =& - \rho_1 \frac{Q^{\text{shell}}}{4\pi}\rmi k_{1} (2n + 1) h_{n} (k_{1} r^{\text{shell}}_{s}) j_n(k_{1} a_{\text{core}}) + \rho_2 \frac{Q^{\text{core}}}{4\pi}\rmi k_{2} (2n + 1) h_{n} (k_{2} a_{\text{core}}) j_n(k_{2} r^{\text{core}}_{s}) \\
        & C_{n} h_{n}'(k_{0} a_{\text{shell}}) - [D_{n} j_{n}'(k_{1} a_{\text{shell}}) + E_{n} y_{n}'(k_{1} a_{\text{shell}})] \nonumber \\
        =& -\frac{Q^{\text{ex}}}{4\pi}\rmi k_{0} (2n + 1) h_{n} (k_{0} r^{\text{ex}}_{s}) j_n'(k_{0} a_{\text{shell}}) + \frac{Q^{\text{shell}}}{4\pi}\rmi k_{1} (2n + 1) h_{n}' (k_{1} a_{\text{shell}}) j_n(k_{1} r^{\text{shell}}_{s}) \\
        & D_{n} j_{n}'(k_{1} a_{\text{core}}) + E_{n} y_{n}'(k_{1} a_{\text{core}}) - F_{n} j_{n}'(k_{2} a_{\text{core}}) \nonumber \\
        =& - \frac{Q^{\text{shell}}}{4\pi}\rmi k_{1} (2n + 1) h_{n} (k_{1} r^{\text{shell}}_{s}) j_n'(k_{1} a_{\text{core}}) + \frac{Q^{\text{core}}}{4\pi}\rmi k_{2} (2n + 1) h_{n}' (k_{2} a_{\text{core}}) j_n(k_{2} r^{\text{core}}_{s}).
    \end{align}
\end{subequations}
when $n=0,1,..., N$. Using the expressions in Eq.~(\ref{eq:anasBf012}) and the recurrence relations in Eq.~(\ref{eq:anasBf}), we can solve Eq.~(\ref{eq:anabcCDEF}) to get coefficients $C_{n} \, D_{n}, \, E_{n}, \, F_{n}$, and we then are able to obtain the acoustic potential field driven by the acoustic monopoles located along the symmetric axis of a concentric core-shell spherical scatterer using the expressions in Eqs.~(\ref{eq:anaex}) to (\ref{eq:anacore}). 

\bibliographystyle{siam} 
\bibliography{sampbib}

\end{document}